\documentclass[11pt,a4paper]{amsart}
\usepackage[utf8]{inputenc}
\usepackage[T1]{fontenc}
\usepackage{lmodern}
\usepackage{amsmath,amssymb,amsthm}
\usepackage{enumitem}
\usepackage{xcolor}
\usepackage[margin=2.8cm]{geometry}
\usepackage[colorlinks=true,linkcolor=blue!50!black,citecolor=blue!50!black,urlcolor=blue!50!black]{hyperref}

\theoremstyle{plain}
\newtheorem{theorem}{Theorem}[section]

\newtheorem{lemma}[theorem]{Lemma}
\newtheorem{corollary}[theorem]{Corollary}

\newtheorem*{mainthm}{Main Theorem}
\newtheorem*{thmA}{Theorem A}
\newtheorem*{thmB}{Theorem B}
\newtheorem*{thmC}{Theorem C}
\newtheorem*{cormoduli}{Corollary D}
\theoremstyle{definition}
\newtheorem{definition}[theorem]{Definition}
\theoremstyle{remark}
\newtheorem{remark}[theorem]{Remark}

\newcommand{\C}{\mathbb{C}}

\newcommand{\Q}{\mathbb{Q}}
\newcommand{\rl}{\sim_{\mathrm{Lip}}}
\newcommand{\dist}{\operatorname{dist}}

\begin{document}

\title[Bi-Lipschitz right-left equivalence of finite plane-to-plane germs]{Bi-Lipschitz right-left equivalence of finite plane-to-plane map germs admits moduli}

\author{Alexandre Fernandes}
\address{Departamento de Matem\'atica, Universidade Federal do Cear\'a, Fortaleza, CE, Brazil}
\email{alex@mat.ufc.br}

\subjclass[2020]{Primary 32S15, 58K40; Secondary 14B05, 32S05}
\keywords{bi-Lipschitz equivalence, right-left equivalence, map germs, moduli, quasihomogeneous polynomials}

\thanks{The author was partially supported by CNPq-Brazil grant  304700/2021-5 and 405618/2023-9.
}

\date{\today}
\begin{abstract}
We show that the bi-Lipschitz right-left equivalence of finite holomorphic map germs $(\C^2,0)\to(\C^2,0)$ admits continuous moduli, already among finitely determined polynomial germs of corank one and degree five. The moduli come from a new invariant of finite germs: for two tangent branches of the critical set whose images are tangent, the limit of the ratio of the norms of the images of nearby points is preserved by bi-Lipschitz right-left equivalence. For germs $(F,H)$ with quasihomogeneous components the invariant is computed from the height functions of $F$ and $H$: it is given by the moduli of the values of the height function of the component of lower weighted degree at the roots of a weighted Wronskian. As a first application we obtain a continuous invariant for the germs $(Y^3+3tYX^4+X^6,\,X^N)$, $N\ge7$, built from the Henry--Parusiński family.
\end{abstract}

\maketitle

\section{Introduction}

For function germs $(\C^2,0)\to(\C,0)$ under bi-Lipschitz \emph{right} equivalence, Henry and Parusiński \cite{HP} discovered that continuous moduli exist; their invariant is built from the leading coefficients of the expansions of the function along the branches of its polar curve. For quasihomogeneous, non-homogeneous polynomials in two variables, bi-Lipschitz right equivalence even coincides with analytic equivalence \cite{CF} (see also \cite{FR,CR}). On the other hand, bi-Lipschitz \emph{contact} equivalence has no moduli \cite{BCFR,RV}, and for polynomial maps $\C^2\to\C^2$ of bounded degree there are only finitely many \emph{topological} right-left types \cite{AN}.

In this paper we consider the bi-Lipschitz right-left equivalence of finite map germs from the plane to the plane.

\begin{definition}
Two germs $P,Q\colon(\C^2,0)\to(\C^2,0)$ are \emph{bi-Lipschitz right-left equivalent}, written $P\rl Q$, if there exist germs of bi-Lipschitz homeomorphisms $\Phi,\Psi\colon(\C^2,0)\to(\C^2,0)$ such that $\Psi\circ P=Q\circ\Phi$.
\end{definition}

Sufficient conditions for this equivalence are given in \cite{CNR,CSS}, and the lowest-degree homogeneous parts of equivalent germs are equivalent \cite{Ng}; these results do not provide continuous invariants. For function germs $(\C^n,0)\to(\C,0)$, Nguyen and Trivedi \cite{NT} showed that the simple germs for the bi-Lipschitz right-left equivalence are the same as for the right equivalence \cite{NRT}; in particular, the family of \cite{HP} is not simple for the right-left equivalence either.

Through the germs $(f,0)$ this gives information on non-finite germs $(\C^2,0)\to(\C^2,0)$ (Remark~\ref{rem:degenerate}), but not on finite ones. For finite germs, to our knowledge, nothing was known about the existence of moduli. The homeomorphism $\Psi$ in the target destroys the Henry--Parusiński invariant in many cases (Remark~\ref{rem:role}), so it is not clear a priori whether moduli survive. We show that they do.

\begin{remark}\label{rem:degenerate}
For some non-finite germs the problem reduces to function germs. Let $f$ and $g$ be non-constant holomorphic germs $(\C^2,0)\to(\C,0)$. Then $(f,0)\rl(g,0)$ if and only if $f$ and $g$ are bi-Lipschitz right-left equivalent as function germs, that is, $\psi\circ f=g\circ\Phi$ for germs of bi-Lipschitz homeomorphisms $\Phi$ of $(\C^2,0)$ and $\psi$ of $(\C,0)$. Indeed, if $\Psi\circ(f,0)=(g,0)\circ\Phi$ then, $f$ and $g$ being open, $\Psi$ maps the germ of the axis $\C\times\{0\}$ onto itself, and $\psi(u)=\Psi_1(u,0)$ satisfies $\psi\circ f=g\circ\Phi$; conversely, one takes $\Psi(u,v)=(\psi(u),v)$. Hence moduli for the bi-Lipschitz right-left equivalence of function germs give moduli among the germs $(f,0)$. No such reduction is available for finite germs, which are the subject of this paper.
\end{remark}

\subsection{The invariant}
A holomorphic germ $P\colon(\C^2,0)\to(\C^2,0)$ is \emph{finite} if $P^{-1}(0)=\{0\}$. We denote by $\Sigma(P)$ its critical set, the zero set of the Jacobian determinant, and call \emph{branches} of $\Sigma(P)$ its irreducible components at $0$. Two irreducible curve germs are \emph{tangent} if they have the same tangent line at $0$. If $\gamma$ is a branch of $\Sigma(P)$, then $P(\gamma)$ is an irreducible curve germ in the target.

\begin{definition}\label{def:Lambda}
Let $P$ be a finite germ. A pair $(\gamma_1,\gamma_2)$ of branches of $\Sigma(P)$ is \emph{admissible} if $\gamma_1,\gamma_2$ are tangent and the curves $P(\gamma_1),P(\gamma_2)$ are tangent. For an admissible pair we set
\[
\Lambda_P(\gamma_1,\gamma_2)=\lim\frac{|P(p_1)|}{|P(p_2)|}\ \in[0,\infty],
\]
the limit being taken over points $p_i\in\gamma_i\setminus\{0\}$ with $p_i\to0$ and $|p_1-p_2|=o(|p_1|)$.
\end{definition}

The limit exists and does not depend on the choice of the points; see \eqref{eq:Lambda} in Section~\ref{sec:prelim}.

\begin{mainthm}
Let $P,Q\colon(\C^2,0)\to(\C^2,0)$ be finite holomorphic germs and $\Phi,\Psi$ germs of bi-Lipschitz homeomorphisms with $\Psi\circ P=Q\circ\Phi$. Then $\Phi$ maps each branch of $\Sigma(P)$ onto a branch of $\Sigma(Q)$, it maps admissible pairs to admissible pairs, and
\[
\Lambda_Q\big(\Phi(\gamma_1),\Phi(\gamma_2)\big)=\Lambda_P(\gamma_1,\gamma_2)
\]
for every admissible pair $(\gamma_1,\gamma_2)$ of branches of $\Sigma(P)$.
\end{mainthm}

Neither the critical set nor the discriminant carries moduli as an embedded curve: the ambient bi-Lipschitz type of a plane curve germ is determined by its embedded topology \cite{NP}; see also \cite{PT,F}.
The number $\Lambda$ is not an invariant of these curves; it measures how $P$ couples the metric of the source with that of the target along them. The proof is elementary. It uses the topological invariance of the critical set of a finite holomorphic germ, the fact that tangency of curves is a bi-Lipschitz notion, and a scaling argument for the restriction of $\Psi$ to the common tangent line of the two image curves.

\subsection{Quasihomogeneous pairs}
We follow the conventions of \cite{CF}. Fix an integer $q\ge2$ and the weights $(1,q)$. A polynomial $F\in\C[X,Y]$ is quasihomogeneous of weighted degree $\nu=\nu(F)$ if $F(\lambda X,\lambda^qY)=\lambda^\nu F(X,Y)$; its \emph{height function} is $f(z)=F(1,z)$, and
\begin{equation}\label{eq:height}
F(X,Y)=X^{\nu}f(Y/X^q)\ \ (X\ne0),\qquad F(s,s^qz)=s^\nu f(z).
\end{equation}

\begin{definition}\label{def:pair}
A \emph{quasihomogeneous pair} with weights $(1,q)$ is a pair $(F,H)$ of quasihomogeneous polynomials with $\nu(F)<\nu(H)$ such that the germ $P=(F,H)\colon(\C^2,0)\to(\C^2,0)$ is finite. The polynomial $F$ is the \emph{leading component} of the pair and its height function $f$ is the \emph{leading height}. The pair is \emph{non-degenerate} if its leading height is not constant, that is, if $F$ is not a constant multiple of $X^{\nu(F)}$.
\end{definition}

For a quasihomogeneous pair $(F,H)$ we always write $P=(F,H)$, $\nu=\nu(F)$, $m=\nu(H)$, $f,h$ for the height functions, and we define the \emph{weighted Wronskian}
\[
W_{P}(z)=\nu\,f(z)h'(z)-m\,f'(z)h(z).
\]
We shall see (Lemma~\ref{lem:sigma}) that $W_P\not\equiv0$ and that, away from $\{X=0\}$, the critical set of $P$ is the union of the curves $\Gamma_i=\{Y=\kappa_iX^q\}$, where $\kappa_1,\dots,\kappa_\ell$ are the distinct roots of $W_P$. We put $v_i=f(\kappa_i)$. For a second pair $(G,K)$, with $Q=(G,K)$, we write $\nu'=\nu(G)$, $m'=\nu(K)$, $g,k$ for the height functions, $\tau_1,\dots,\tau_{\ell'}$ for the distinct roots of $W_Q$, and $w_j=g(\tau_j)$.

When $H=X^N$ we have $h\equiv1$ and $W_P=-Nf'$: the $\kappa_i$ are the critical points of $f$ and the $v_i$ its critical values, which are the Henry--Parusiński invariants of $F$ \cite[Corollary 7.9]{CF}.

\begin{thmA}
Let $(F,H)$ and $(G,K)$ be non-degenerate quasihomogeneous pairs with the same weights $(1,q)$, and assume that $W_P$ has at least two distinct roots. If $P=(F,H)\rl Q=(G,K)$, then:
\begin{enumerate}[label=\textup{(\alph*)}]
\item $\nu(F)=\nu(G)$ and $\nu(H)=\nu(K)$.
\item $\ell=\ell'$ and there is a bijection $\sigma$ of $\{1,\dots,\ell\}$ such that $v_i=0\iff w_{\sigma(i)}=0$.
\item There is $\rho>0$ such that $|w_{\sigma(i)}|=\rho\,|v_i|$ for all $i$.
\end{enumerate}
\end{thmA}

Theorem~A follows from the Main Theorem: for $v_i,v_j\ne0$ the pair $(\Gamma_i,\Gamma_j)$ is admissible and $\Lambda_P(\Gamma_i,\Gamma_j)=|v_i|/|v_j|$. The condition $\nu(F)<\nu(H)$ in Definition~\ref{def:pair} is not a restriction on the order of the components, since one may exchange the coordinates of the target; what matters is non-degeneracy, that is, that the height function of the component of \emph{lower} degree is not constant (Remark~\ref{rem:role}).

\subsection{Germs of the form $(F,X^N)$}
Let
\[
F_t(X,Y)=Y^3+3tYX^4+X^6,\qquad f_t(z)=z^3+3tz+1,\qquad P_t=(F_t,X^N),
\]
with weights $(1,2)$ and weighted degree $6$; this is the family $x^3-3\tau^2xy^4+y^6$ of \cite{HP}, with $t=-\tau^2$. The critical points of $f_t$ are $\kappa_\pm=\pm(-t)^{1/2}$ and the critical values are
\[
v_\pm(t)=1\pm\theta_t,\qquad \theta_t=2t\,(-t)^{1/2},\qquad \theta_t^2=-4t^3.
\]
Let $T=\{t\in\C:\ t\ne0,\ 4t^3\ne-1\}$ and, for $t\in T$,
\[
\lambda_t=\frac{1+\theta_t}{1-\theta_t},\qquad I(t)=\big|\log|\lambda_t|\big|\in[0,\infty).
\]
Changing the sign of $\theta_t$ replaces $\lambda_t$ by $\lambda_t^{-1}$, so $I(t)$ depends only on $t^3$.

\begin{thmB}
Let $N\ge7$ and $t,s\in T$. If $P_t\rl P_s$ then $I(t)=I(s)$.
\end{thmB}

The function $I$ takes all the values in $[0,\infty)$. For comparison, $F_t$ and $F_s$ are bi-Lipschitz right equivalent if and only if $t^3=s^3$ \cite{HP,CF}. The restriction $N\ge7$ is essential: for $N\in\{1,2,3,6\}$ all the germs $P_t$, $t\ne0$, are analytically right-left equivalent (Remark~\ref{rem:role}).

\subsection{A finitely determined family of corank one}
For $\lambda\in\C\setminus\{0,1\}$ let
\[
H_\lambda(X,Y)=15XY^2-5(1+\lambda)X^3Y+3\lambda X^5,\qquad P_\lambda=(Y,\,H_\lambda),
\]
with weights $(1,2)$, $\nu(Y)=2<5=\nu(H_\lambda)$. Its critical set is $\{\partial H_\lambda/\partial X=0\}$, and $\partial H_\lambda/\partial X=15(Y-X^2)(Y-\lambda X^2)$: two tangent parabolas.

\begin{thmC}
Let $\lambda,\lambda'\in\C\setminus\{0,1\}$.
\begin{enumerate}[label=\textup{(\roman*)}]
\item If $P_\lambda\rl P_{\lambda'}$ then $\big|\log|\lambda'|\big|=\big|\log|\lambda|\big|$.
\item For $\lambda$ outside a finite set, $P_\lambda$ is finitely determined \textup{(}for analytic right-left equivalence\textup{)}.
\end{enumerate}
\end{thmC}

\begin{cormoduli}
The bi-Lipschitz right-left equivalence of finite holomorphic germs $(\C^2,0)\to(\C^2,0)$ admits continuous moduli. More precisely, there are uncountably many classes among polynomial germs of degree $5$ which are finitely determined and of corank one, and among the finite polynomial germs $(F_t,X^7)$ of degree $7$.
\end{cormoduli}

All these germs belong to finitely many topological right-left classes \cite{AN}.

\subsection{Organization}
Section~\ref{sec:prelim} contains the preliminaries. The Main Theorem is proved in Section~\ref{sec:main}. Section~\ref{sec:qh} describes the critical set of a quasihomogeneous pair and proves Theorem~A. Theorems B and C are proved in Section~\ref{sec:ex}.

\section{Preliminaries}\label{sec:prelim}

The results of this section are classical. We state them in the form in which they will be used and refer to the literature for the proofs.

\subsection{Multiplicity}
Let $P$ be a holomorphic map between open subsets of $\C^2$ which is \emph{finite at $p$}, that is, $p$ is an isolated point of its fibre. The \emph{multiplicity} $\mu_p(P)$ of $P$ at $p$ is the dimension of its local algebra at $p$. It coincides with the local topological degree of $P$ at $p$: for every small neighbourhood $U$ of $p$, each point close to $P(p)$ has at most $\mu_p(P)$ preimages in $U$, and a generic one has exactly $\mu_p(P)$; see \cite[Chapters 4--5]{AGV}. Moreover $\mu_p(P)=1$ if and only if $p\notin\Sigma(P)$.

Since homeomorphisms have local degree $\pm1$, the multiplicity is preserved by topological right-left equivalence, and so are the critical set and the discriminant; see for instance the survey \cite{EH} and the references therein. For the critical set this is also a consequence of Osgood's theorem \cite{Na}, by which a holomorphic map between domains of $\C^2$ is locally injective exactly at the points where its Jacobian does not vanish.

\begin{lemma}[see {\cite{AGV,EH}}]\label{lem:top}
Let $P,Q\colon(\C^2,0)\to(\C^2,0)$ be finite holomorphic germs and $\Phi,\Psi$ germs of homeomorphisms with $\Psi\circ P=Q\circ\Phi$. Then $\mu_{\Phi(p)}(Q)=\mu_p(P)$ for every $p$ close to $0$. In particular $\Phi(\Sigma(P))=\Sigma(Q)$.
\end{lemma}

\subsection{Branches and tangency}
Let $\gamma$ be an irreducible curve germ at $0\in\C^2$. By Puiseux's theorem it has an injective holomorphic parametrization $p(\tau)=a\,\tau^n+O(\tau^{n+1})$, $a\ne0$, and the line $T(\gamma)=\C a$ is its tangent line (its tangent cone) at $0$. The next lemma collects two standard consequences of the Puiseux parametrization; see for instance \cite[Chapters 2 and 4]{Wallbook}, and \cite{FS} for the metric point of view adopted here.

\begin{lemma}[see {\cite{Wallbook,FS}}]\label{lem:branch}
Let $\gamma,\gamma'$ be irreducible curve germs at $0\in\C^2$.
\begin{enumerate}[label=\textup{(\arabic*)}]
\item $\dist(p,T(\gamma))=o(|p|)$ for $p\in\gamma$, and the orthogonal projection of $\gamma$ to $T(\gamma)$ covers a neighbourhood of $0$ in $T(\gamma)$.
\item $\gamma$ and $\gamma'$ are tangent if and only if $\dist(p,\gamma')=o(|p|)$ for $p\in\gamma$.
\end{enumerate}
\end{lemma}

In other words, two distinct branches are tangent if and only if their order of contact is greater than $1$. The order of contact of pairs of branches is preserved by bi-Lipschitz homeomorphisms \cite{PT,F,NP} (see also \cite{FS}). We only need the following particular case, which also follows from the bi-Lipschitz invariance of tangent cones \cite{S}.

\begin{lemma}[tangency is a bi-Lipschitz notion; {\cite{PT,F,NP,S}}]\label{lem:tangent}
Let $\gamma_1,\gamma_2$ be irreducible curve germs at $0\in\C^2$ and let $\Theta$ be a germ of bi-Lipschitz homeomorphism of $(\C^2,0)$ such that $\Theta(\gamma_1),\Theta(\gamma_2)$ are irreducible curve germs. Then $\gamma_1,\gamma_2$ are tangent if and only if $\Theta(\gamma_1),\Theta(\gamma_2)$ are tangent.
\end{lemma}

\subsection{Orders along a branch}
Let $P\colon(\C^2,0)\to(\C^2,0)$ be a holomorphic germ which is not identically zero on $\gamma$, and write $P(p(\tau))=b\,\tau^k+O(\tau^{k+1})$ with $b\ne0$. Then $P(\gamma)$ is an irreducible curve germ with tangent line $\C b$, and
\begin{equation}\label{eq:order}
|P(p)|=C\,|p|^{e}\,(1+o(1))\qquad(p\in\gamma,\ p\to0),
\end{equation}
where $e=e_P(\gamma)=k/n\in\Q_{>0}$ and $C=C_P(\gamma)=|b|/|a|^{e}>0$. Now let $P$ be finite and let $\gamma_1,\gamma_2$ be tangent irreducible curve germs, with $e_i=e_P(\gamma_i)$ and $C_i=C_P(\gamma_i)$. If $p_i\in\gamma_i\setminus\{0\}$ are points with $p_i\to0$ and $|p_1-p_2|=o(|p_1|)$, then $|p_2|/|p_1|\to1$ and \eqref{eq:order} gives
\begin{equation}\label{eq:Lambda}
\frac{|P(p_1)|}{|P(p_2)|}\longrightarrow\Lambda_P(\gamma_1,\gamma_2)=
\begin{cases}
C_1/C_2&\text{if }e_1=e_2,\\
0&\text{if }e_1>e_2,\\
\infty&\text{if }e_1<e_2.
\end{cases}
\end{equation}
In particular, the limit in Definition~\ref{def:Lambda} exists and does not depend on the choice of the points.

\section{Proof of the Main Theorem}\label{sec:main}

The following lemma is an approximate version of a simple remark: if a positive function $\mathsf n$ on $\C\setminus\{0\}$ is comparable to $|z|$ and satisfies $\mathsf n(\Lambda z)=\Lambda'\,\mathsf n(z)$, then $\Lambda=\Lambda'$, because iterating the relation $k$ times multiplies the bounded quotient $\mathsf n(z)/|z|$ by $(\Lambda'/\Lambda)^k$.

\begin{lemma}[scaling lemma]\label{lem:scaling}
Let $\mathsf n$ be a positive function on a punctured disc $\{0<|z|<r_0\}\subset\C$ with $c^{-1}|z|\le\mathsf n(z)\le c\,|z|$ for some $c\ge1$, and let $\Lambda,\Lambda'\in(0,\infty)$. Suppose that for every $\eta>0$ there is $\delta>0$ such that for every $z$ with $0<|z|<\delta$ there is $\tilde z$ with
\[
\Big|\frac{|\tilde z|}{\Lambda|z|}-1\Big|<\eta\qquad\text{and}\qquad\Big|\frac{\mathsf n(\tilde z)}{\Lambda'\,\mathsf n(z)}-1\Big|<\eta .
\]
Then $\Lambda=\Lambda'$.
\end{lemma}

\begin{proof}
Let $\varrho=\Lambda'/\Lambda$ and suppose $\varrho\ne1$; put $\varrho_*=\max\{\varrho,\varrho^{-1}\}>1$. Choose $\eta\in(0,1)$ with $(1+\eta)/(1-\eta)<\varrho_*^{1/2}$, then $n$ with $\varrho_*^{n/2}>c^2$, then $\delta$ as in the hypothesis, and $z_0$ with $0<|z_0|<\delta\min\{1,(\Lambda(1+\eta))^{-n}\}$. Define $z_{k+1}=\tilde z_k$ for $k<n$; then $|z_k|\le|z_0|(\Lambda(1+\eta))^k<\delta$, so the construction is possible, and
\[
\frac{|z_n|}{|z_0|}\in\Lambda^n\big[(1-\eta)^n,(1+\eta)^n\big],\qquad \frac{\mathsf n(z_n)}{\mathsf n(z_0)}\in\Lambda'^n\big[(1-\eta)^n,(1+\eta)^n\big].
\]
The quotient $\Xi=\big(\mathsf n(z_n)/|z_n|\big)\big/\big(\mathsf n(z_0)/|z_0|\big)$ lies in $[c^{-2},c^2]$. If $\varrho>1$ then $\Xi\ge\varrho^n\big((1-\eta)/(1+\eta)\big)^n>\varrho^{n/2}>c^2$; if $\varrho<1$ then $\Xi\le\varrho^n\big((1+\eta)/(1-\eta)\big)^n<\varrho^{n/2}<c^{-2}$. Both are impossible.
\end{proof}

\begin{proof}[Proof of the Main Theorem]
Let $c\ge1$ be such that $\Phi,\Psi$ and their inverses are $c$-Lipschitz near $0$.

\emph{Step 1: branches.} By Lemma~\ref{lem:top}, $\Phi(\Sigma(P))=\Sigma(Q)$. Near $0$, the set $\Sigma(P)\setminus\{0\}$ is the disjoint union of the sets $\gamma\setminus\{0\}$, $\gamma$ a branch; each of them is connected (it is the image of a punctured disc by the parametrization) and open and closed in $\Sigma(P)\setminus\{0\}$. Hence $\Phi$ maps $\gamma\setminus\{0\}$ into $\gamma'\setminus\{0\}$ for a single branch $\gamma'$ of $\Sigma(Q)$; applying the same to $\Phi^{-1}$ we get $\Phi(\gamma)=\gamma'$ as germs, and $\gamma\mapsto\gamma'$ is a bijection between the sets of branches.

\emph{Step 2: admissible pairs.} Let $\gamma_1,\gamma_2$ be branches of $\Sigma(P)$ and $\gamma'_i=\Phi(\gamma_i)$. By Lemma~\ref{lem:tangent} applied to $\Phi$, $\gamma_1,\gamma_2$ are tangent if and only if $\gamma'_1,\gamma'_2$ are. Since $\Psi(P(\gamma_i))=Q(\gamma'_i)$, the same lemma applied to $\Psi$ shows that $P(\gamma_1),P(\gamma_2)$ are tangent if and only if $Q(\gamma'_1),Q(\gamma'_2)$ are.

\emph{Step 3: orders.} For $p\in\gamma_i$ we have $|Q(\Phi(p))|=|\Psi(P(p))|\in[c^{-1},c]\,|P(p)|$ and $|\Phi(p)|\in[c^{-1},c]\,|p|$. By \eqref{eq:order}, $|P(p)|$ is comparable to $|p|^{e_P(\gamma_i)}$ and $|Q(\Phi(p))|$ to $|\Phi(p)|^{e_Q(\gamma'_i)}$; hence $e_Q(\gamma'_i)=e_P(\gamma_i)$. By \eqref{eq:Lambda} this proves the theorem when $e_P(\gamma_1)\ne e_P(\gamma_2)$.

\emph{Step 4: the case of equal orders.} Let $(\gamma_1,\gamma_2)$ be admissible with $e_P(\gamma_1)=e_P(\gamma_2)$, and put $\Lambda=\Lambda_P(\gamma_1,\gamma_2)$, $\Lambda'=\Lambda_Q(\gamma'_1,\gamma'_2)$, both in $(0,\infty)$. Let $D_i=P(\gamma_i)$, let $L$ be the common tangent line of $D_1,D_2$, $\epsilon\in L$ a unit vector, and for $y\in\C^2$ write $\zeta(y)\epsilon$ for the orthogonal projection of $y$ to $L$. Define
\[
\mathsf n(z)=|\Psi(z\epsilon)|,\qquad c^{-1}|z|\le\mathsf n(z)\le c\,|z| .
\]
For $y\in D_1\cup D_2$ we have $|y-\zeta(y)\epsilon|=o(|y|)$ by Lemma~\ref{lem:branch}(1); hence $|\zeta(y)|=|y|(1+o(1))$ and, since $\Psi$ is $c$-Lipschitz and $\mathsf n(\zeta(y))\ge|\zeta(y)|/c$,
\begin{equation}\label{eq:npsi}
|\Psi(y)|=\mathsf n(\zeta(y))\,(1+o(1))\qquad(y\in D_1\cup D_2,\ y\to0).
\end{equation}
We verify the hypothesis of Lemma~\ref{lem:scaling}. Let $z\ne0$ be small. By Lemma~\ref{lem:branch}(1) applied to $D_2$ there is $p_2\in\gamma_2$ with $\zeta(P(p_2))=z$, and $p_2\to0$ as $z\to0$. Since $\gamma_1,\gamma_2$ are tangent, by Lemma~\ref{lem:branch}(2) there is $p_1\in\gamma_1$ with $|p_1-p_2|=o(|p_2|)$. Put $\tilde z=\zeta(P(p_1))$. Then, by \eqref{eq:Lambda},
\[
\frac{|\tilde z|}{|z|}=\frac{|P(p_1)|}{|P(p_2)|}\,(1+o(1))\longrightarrow\Lambda .
\]
On the other hand, by \eqref{eq:npsi} and $\Psi\circ P=Q\circ\Phi$,
\[
\frac{\mathsf n(\tilde z)}{\mathsf n(z)}=\frac{|\Psi(P(p_1))|}{|\Psi(P(p_2))|}\,(1+o(1))=\frac{|Q(\Phi(p_1))|}{|Q(\Phi(p_2))|}\,(1+o(1))\longrightarrow\Lambda' ,
\]
because $\Phi(p_i)\in\gamma'_i$ and $|\Phi(p_1)-\Phi(p_2)|\le c\,|p_1-p_2|=o(|\Phi(p_2)|)$, so that \eqref{eq:Lambda} applies to $Q$. All the $o(1)$ are uniform as $z\to0$. By Lemma~\ref{lem:scaling}, $\Lambda=\Lambda'$.
\end{proof}

\begin{remark}\label{rem:preimage}
The proof uses only that $\Phi$ maps the curves $\gamma_i$ onto analytic branches. Since $\Phi$ also maps $P^{-1}(P(\Sigma(P)))$ onto $Q^{-1}(Q(\Sigma(Q)))$, the Main Theorem holds, with the same proof, for admissible pairs of branches of the preimage of the discriminant.
\end{remark}

\section{Quasihomogeneous pairs}\label{sec:qh}

\begin{lemma}[critical set]\label{lem:sigma}
Let $(F,H)$ be a non-degenerate quasihomogeneous pair, $\nu=\nu(F)$, $m=\nu(H)$.
\begin{enumerate}[label=\textup{(\arabic*)}]
\item For $s\ne0$ and $z\in\C$, the Jacobian $J_P=\partial_XF\,\partial_YH-\partial_YF\,\partial_XH$ satisfies
\[
J_P(s,s^qz)=s^{\nu+m-q-1}\,W_P(z).
\]
\item $W_P\not\equiv0$, and $(f(\kappa),h(\kappa))\ne(0,0)$ for every $\kappa\in\C$.
\item $\Sigma(P)\cap\{X\ne0\}=\bigcup_{i=1}^{\ell}\Gamma_i\setminus\{0\}$, where $\Gamma_i=\{Y=\kappa_iX^q\}$. Thus the branches of $\Sigma(P)$ are $\Gamma_1,\dots,\Gamma_\ell$ and possibly $L=\{X=0\}$.
\item For $p=(s,\kappa_is^q)\in\Gamma_i$ we have $|p|=|s|(1+o(1))$ and $P(p)=(s^\nu v_i,\,s^mh(\kappa_i))$. Hence, if $v_i\ne0$, then $e_P(\Gamma_i)=\nu$, $C_P(\Gamma_i)=|v_i|$ and $P(\Gamma_i)$ is tangent to the axis $\{v=0\}$; if $v_i=0$, then $e_P(\Gamma_i)=m$ and $P(\Gamma_i)$ is contained in the axis $\{u=0\}$.
\end{enumerate}
\end{lemma}

\begin{proof}
(1) Differentiating $F(X,Y)=X^\nu f(Y/X^q)$ on $X\ne0$, with $z=Y/X^q$,
\[
\partial_YF=X^{\nu-q}f'(z),\qquad \partial_XF=\nu X^{\nu-1}f(z)-qYX^{\nu-q-1}f'(z)=X^{\nu-1}\big(\nu f(z)-qzf'(z)\big),
\]
and similarly $\partial_YH=X^{m-q}h'(z)$, $\partial_XH=X^{m-1}(mh(z)-qzh'(z))$. Therefore
\[
J_P=X^{\nu+m-q-1}\Big[(\nu f-qzf')h'-f'(mh-qzh')\Big]=X^{\nu+m-q-1}\big(\nu fh'-mf'h\big).
\]

(2) If $f(\kappa)=h(\kappa)=0$, then $P(s,s^q\kappa)=(s^\nu f(\kappa),s^mh(\kappa))=0$ for all $s$, and $P$ is not finite. If $h\equiv0$ then $H\equiv0$, and again $P$ is not finite. Suppose $W_P\equiv0$. On the open set where $fh\ne0$ we have $(h^\nu/f^m)'=h^{\nu-1}f^{-m-1}W_P=0$, hence $h^\nu=c\,f^m$ for a constant $c$, an identity valid on $\C$. Since $f$ is not constant, it has a root $\kappa$, and then $h(\kappa)=0$: a common root, which is impossible.

(3) Every point with $X\ne0$ is uniquely of the form $(s,s^qz)$, and by (1) it belongs to $\Sigma(P)$ if and only if $W_P(z)=0$. On $L$: the polynomial $J_P$ is quasihomogeneous, so $J_P(0,Y)$ is a monomial $c\,Y^{e}$; if $c=0$ then $L\subset\Sigma(P)$, and if $c\ne0$ then $\Sigma(P)\cap L=\{0\}$.

(4) is immediate from \eqref{eq:height}, (2) and $\nu<m$.
\end{proof}

\begin{lemma}[weighted Bézout]\label{lem:bezout}
If $F,H$ are quasihomogeneous with weights $(1,q)$ and $P=(F,H)$ is finite, then $\mu_0(P)=\nu(F)\,\nu(H)/q$.
\end{lemma}

This is the classical Bézout formula for quasihomogeneous maps; see \cite[Chapter~12]{AGV}.

We include a short proof for the convenience of the reader.

\begin{proof}
Recall that $\mu_0(P)$ is the intersection multiplicity $i_0(F,H)$ of the curves $\{F=0\}$ and $\{H=0\}$ at the origin. By \eqref{eq:height}, $F=c\,X^{a}\prod_j(Y-\lambda_jX^q)^{m_j}$, where the $\lambda_j$ are the roots of $f$, the $m_j$ their multiplicities, and $\nu(F)=a+q\sum_jm_j$. The intersection multiplicity is additive in $F$. The number $i_0(Y-\lambda_jX^q,H)$ is the order of $H(X,\lambda_jX^q)=X^{\nu(H)}h(\lambda_j)$, which is $\nu(H)$ because $f$ and $h$ have no common root, $P$ being finite. If $a>0$, then $X\nmid H$, so $H(0,Y)$ is a non-zero multiple of $Y^{\nu(H)/q}$ and $i_0(X,H)=\nu(H)/q$. Hence
\[
i_0(F,H)=a\,\frac{\nu(H)}{q}+\nu(H)\sum_jm_j=\frac{\nu(F)\,\nu(H)}{q}.\qedhere
\]
\end{proof}

\begin{corollary}[correspondence of branches]\label{cor:branches}
Under the hypotheses of Theorem~A, $W_P$ and $W_Q$ have the same number $\ell$ of distinct roots, and there is a bijection $\sigma$ of $\{1,\dots,\ell\}$ such that $\Phi(\Gamma_i)=\Gamma'_{\sigma(i)}$ as germs, where $\Gamma'_j=\{Y=\tau_jX^q\}$.
\end{corollary}

\begin{proof}
By Step 1 of the proof of the Main Theorem, $\Phi$ induces a bijection between the branches of $\Sigma(P)$ and those of $\Sigma(Q)$, which by Lemma~\ref{lem:sigma}(3) are the $\Gamma_i$ and possibly $L$, respectively the $\Gamma'_j$ and possibly $L'=\{X=0\}$. Let $i\ne j$. The curves $\Gamma_i,\Gamma_j$ are tangent (their tangent line is $\{Y=0\}$, as $q\ge2$), so by Lemma~\ref{lem:tangent} $\Phi(\Gamma_i)$ and $\Phi(\Gamma_j)$ are distinct tangent branches of $\Sigma(Q)$. But $L'$ is not tangent to any other branch of $\Sigma(Q)$, all of which have tangent line $\{Y=0\}$. Hence $\Phi(\Gamma_i)\ne L'$ for all $i$. In particular $W_Q$ has at least two distinct roots, and the same argument applies to $\Phi^{-1}$.
\end{proof}

\begin{proof}[Proof of Theorem A]
Let $\sigma$ be the bijection of Corollary~\ref{cor:branches}.

By Lemma~\ref{lem:top} at $p=0$ and Lemma~\ref{lem:bezout},
\begin{equation}\label{eq:bez}
\nu m=\nu'm' .
\end{equation}
By Step 3 of the proof of the Main Theorem, $e_Q(\Gamma'_{\sigma(i)})=e_P(\Gamma_i)$ for all $i$, and by Lemma~\ref{lem:sigma}(4) these orders are $\nu$ or $m$ (resp.\ $\nu'$ or $m'$) according to whether $v_i\ne0$ or $v_i=0$ (resp.\ $w_{\sigma(i)}\ne0$ or $w_{\sigma(i)}=0$). Suppose that $v_i\ne0$ and $w_{\sigma(i)}=0$ for some $i$. Then $\nu=m'$, and \eqref{eq:bez} gives $m=\nu'$; but then $m=\nu'<m'=\nu<m$, a contradiction. The case $v_i=0$, $w_{\sigma(i)}\ne0$ is excluded in the same way. Hence $v_i=0\iff w_{\sigma(i)}=0$. If some $v_i\ne0$ we get $\nu=\nu'$ from the orders and then $m=m'$ from \eqref{eq:bez}; if all $v_i$ vanish we get $m=m'$ and then $\nu=\nu'$. This proves (a) and (b).

(c) Let $i,j$ be such that $v_i,v_j\ne0$. By Lemma~\ref{lem:sigma}(4) the pair $(\Gamma_i,\Gamma_j)$ is admissible and $\Lambda_P(\Gamma_i,\Gamma_j)=|v_i|/|v_j|$; likewise $\Lambda_Q(\Gamma'_{\sigma(i)},\Gamma'_{\sigma(j)})=|w_{\sigma(i)}|/|w_{\sigma(j)}|$. By the Main Theorem the two numbers coincide. Fixing $j_0$ with $v_{j_0}\ne0$ and $\rho=|w_{\sigma(j_0)}|/|v_{j_0}|$ we get $|w_{\sigma(i)}|=\rho|v_i|$ for all $i$; if all $v_i$ vanish there is nothing to prove.
\end{proof}

\begin{remark}[which coordinate carries the invariant]\label{rem:role}
In Theorem~A the invariant is read on the component of \emph{lower} weighted degree. If $\nu(H)<\nu(F)$ one exchanges the coordinates of the target and applies the theorem to $(H,F)$. If the height function of the component of lower degree is constant, the theorem gives nothing. This is the case of $(F,X^N)$ with $N<\nu(F)$, that is, of the pair $(X^N,F)$. It is consistent with what happens when $N$ divides $\nu$: the change $(u,v)\mapsto(u-bv^{\nu/N},v)$ in the target translates all the critical values of $f$. For the Henry--Parusiński family and $N\in\{1,2,3,6\}$ all the germs $(F_t,X^N)$, $t\ne0$, are analytically right-left equivalent: $(u,v)\mapsto(u-v^{6/N},v)$ removes $X^6$, and $Y\mapsto cY$ with $c^2=t$ followed by $u\mapsto u/c^3$ removes $t$. In particular, for unfoldings $(X,F(X,Y))$, whose parameter is the \emph{light} variable, our invariants are trivial, while for unfoldings $(Y,H(X,Y))$, whose parameter is the \emph{heavy} variable, they give Theorem~C.
\end{remark}

\section{The examples}\label{sec:ex}

\begin{proof}[Proof of Theorem B]
$(F_t,X^N)$ is a non-degenerate quasihomogeneous pair with $q=2$ and $\nu=6<N$: it is finite because $F_t(0,Y)=Y^3$, and $f_t$ is not constant. Here $W_{P_t}=-Nf_t'=-3N(z^2+t)$, which for $t\ne0$ has the two simple roots $\kappa_\pm=\pm(-t)^{1/2}$, with
\[
v_\pm(t)=f_t(\kappa_\pm)=\kappa_\pm(\kappa_\pm^2+3t)+1=1+2t\kappa_\pm=1\pm\theta_t .
\]
For $t\in T$ both values are non-zero. If $P_t\rl P_s$ with $t,s\in T$, Theorem~A(c) gives $|v_+(s)/v_-(s)|^{\pm1}=|v_+(t)/v_-(t)|$, the sign depending on the bijection $\sigma$; that is, $|\lambda_s|^{\pm1}=|\lambda_t|$ and $I(s)=I(t)$.
\end{proof}

\begin{remark}
The map $t\mapsto\lambda_t$ takes on $T$ all the values in $\C\setminus\{0,\pm1\}$, so $I$ takes all the values in $[0,\infty)$.
\end{remark}

We turn to the family $P_\lambda=(Y,H_\lambda)$, with $H_\lambda=15XY^2-5(1+\lambda)X^3Y+3\lambda X^5$, weights $(1,2)$, $\nu=2$, $m=5$, and height functions $f(z)=z$, $h_\lambda(z)=15z^2-5(1+\lambda)z+3\lambda$.

\begin{lemma}\label{lem:D}
For $\lambda\in\C\setminus\{0,1\}$, $(Y,H_\lambda)$ is a non-degenerate quasihomogeneous pair and $W_{P_\lambda}=-15(z-1)(z-\lambda)$. Its roots are $\kappa_1=1$, $\kappa_2=\lambda$, with $v_1=1$, $v_2=\lambda$, and $h_\lambda(1)=2(5-\lambda)$, $h_\lambda(\lambda)=2\lambda(5\lambda-1)$.
\end{lemma}

\begin{proof}
$Y$ and $H_\lambda$ are quasihomogeneous with weights $(1,2)$ and degrees $2<5$. If $Y=0$ and $H_\lambda=0$ then $3\lambda X^5=0$, so the germ is finite. The height function $f(z)=z$ is not constant, and
\[
W_{P_\lambda}=2z\,h_\lambda'(z)-5\,h_\lambda(z)=2z\big(30z-5(1+\lambda)\big)-5\,h_\lambda(z)=-15\big(z^2-(1+\lambda)z+\lambda\big).
\]
The values of $h_\lambda$ are obtained by substitution.
\end{proof}

\begin{proof}[Proof of Theorem C]
(i) By Lemma~\ref{lem:D}, Theorem~A applies to $P_\lambda$ and $P_{\lambda'}$. The bijection $\sigma$ maps $\{1,\lambda\}$ to $\{1,\lambda'\}$, and item (c) gives $|\lambda'|^{\pm1}=|\lambda|$.

(ii) Let $\lambda\notin\{0,1,5,1/5\}$ be such that $\lambda^3(5-\lambda)^2\ne(5\lambda-1)^2$. We show that $P_\lambda$ is stable outside the origin; by the Mather--Gaffney geometric criterion \cite{Wall}, this is equivalent to finite determinacy. The Jacobian of $P_\lambda$ is $-\partial H_\lambda/\partial X=-15(Y-X^2)(Y-\lambda X^2)$, so $\Sigma(P_\lambda)=\Gamma_1\cup\Gamma_2$. Moreover $\partial^2H_\lambda/\partial X^2=30X\big(2\lambda X^2-(1+\lambda)Y\big)$ equals $30(\lambda-1)X^3$ on $\Gamma_1$ and $30\lambda(1-\lambda)X^3$ on $\Gamma_2$; since it does not vanish on $\Sigma\setminus\{0\}$, all the points of $\Sigma\setminus\{0\}$ are folds of the map $(Y,H_\lambda)$. The restriction of $P_\lambda$ to $\Gamma_i$ is $s\mapsto(\kappa_is^2,h_\lambda(\kappa_i)s^5)$ with $\kappa_i\ne0$ and $h_\lambda(\kappa_i)\ne0$ (Lemma~\ref{lem:D}), hence injective. The images of $\Gamma_1$ and $\Gamma_2$ are the curves $v^2=\varrho_iu^5$ with $\varrho_i=h_\lambda(\kappa_i)^2/\kappa_i^5$, and $\varrho_1\ne\varrho_2$ is equivalent to $\lambda^3(5-\lambda)^2\ne(5\lambda-1)^2$. Thus, outside the origin, the only singular points are folds and $P_\lambda|_\Sigma$ is injective.
\end{proof}

Note that $H_\lambda(aX,Y)=a\,H_{1/\lambda}(X,Y)$ when $a^2=1/\lambda$, so that $P_\lambda$ and $P_{1/\lambda}$ are analytically right-left equivalent; this is why only $|\log|\lambda||$ appears in Theorem~C.

Corollary~D follows from Theorems B and C, since $|\log|\lambda||$ takes all the values in $[0,\infty)$ even after removing finitely many values of $\lambda$.

\subsection*{Declaration of AI use}

During the development of this work, the author used Claude (Anthropic) as an AI-assisted tool for mathematical discussion, exploration of ideas, and critical examination and refinement of arguments. The author reviewed the resulting mathematical content, independently checked the statements and proofs, and assumes full responsibility for the originality, accuracy, and integrity of the manuscript.


\begin{thebibliography}{99}
\bibitem{AN} K.~Aoki, H.~Noguchi, \emph{On topological types of polynomial map germs of plane to plane}, Mem. School Sci. Engrg. Waseda Univ. \textbf{44} (1980), 133--156.

\bibitem{AGV} V.~I. Arnold, S.~M. Gusein-Zade, A.~N. Varchenko, \emph{Singularities of Differentiable Maps}, vol.~I, Birkhäuser, 1985.

\bibitem{BCFR} L.~Birbrair, J.~C.~F. Costa, A.~Fernandes, M.~A.~S. Ruas, \emph{$\mathcal K$-bi-Lipschitz equivalence of real function-germs}, Proc. Amer. Math. Soc. \textbf{135} (2007), 1089--1095.

\bibitem{CF} L.~Câmara, A.~Fernandes, \emph{Bi-Lipschitz versus analytic equivalence of two variable complex quasihomogeneous function-germs}, preprint, arXiv:2503.02085.

\bibitem{CR} L.~M. Câmara, M.~A.~S. Ruas, \emph{On the moduli space of quasi-homogeneous functions}, Bull. Braz. Math. Soc. (N.S.) \textbf{53} (2022), 895--908.

\bibitem{CNR} J.~C.~F. Costa, T.~Nishimura, M.~A.~S. Ruas, \emph{Bi-Lipschitz $\mathcal A$-equivalence of $\mathcal K$-equivalent map-germs}, Rev. R. Acad. Cienc. Exactas Fís. Nat. Ser. A Mat. RACSAM \textbf{108} (2014), 173--182.

\bibitem{CSS} J.~C.~F. Costa, M.~J. Saia, C.~H. Soares Júnior, \emph{Bi-Lipschitz $\mathcal A$-triviality of map germs and Newton filtrations}, Topology Appl. \textbf{159} (2012), 430--436.

\bibitem{EH} B.~El~Hilany, \emph{Around the topological classification problem of polynomial maps: a survey}, Japanese Journal of Mathematics \textbf{21} (2026), no.~1, 1--66.

\bibitem{F} A.~Fernandes, \emph{Topological equivalence of complex curves and bi-Lipschitz homeomorphisms}, Michigan Math. J. \textbf{51} (2003), 593--606.

\bibitem{FR} A.~Fernandes, M.~A.~S. Ruas, \emph{Rigidity of bi-Lipschitz equivalence of weighted homogeneous function-germs in the plane}, Proc. Amer. Math. Soc. \textbf{141} (2013), 1125--1133.

\bibitem{FS} A.~Fernandes, J.~E. Sampaio, \emph{Lipschitz geometry of complex singularities}, Rio de Janeiro: Editora do IMPA, 2025. 114 p. (Coleção Colóquios Brasileiros de Matemática). ISBN 978-85-244-0592-1.

\bibitem{HP} J.-P. Henry, A.~Parusiński, \emph{Existence of moduli for bi-Lipschitz equivalence of analytic functions}, Compositio Math. \textbf{136} (2003), 217--235.

\bibitem{Na} R.~Narasimhan, \emph{Several Complex Variables}, University of Chicago Press, 1971.

\bibitem{NP} W.~D. Neumann, A.~Pichon, \emph{Lipschitz geometry of complex curves}, J. Singul. \textbf{10} (2014), 225--234.

\bibitem{Ng} N.~Nguyen, \emph{Bi-Lipschitz contact invariance of rank}, Dalat Univ. J. Sci. (2022), 40--47; arXiv:2001.04739.

\bibitem{NRT} N.~Nguyen, M.~Ruas, S.~Trivedi, \emph{Classification of Lipschitz simple function germs}, Proc. Lond. Math. Soc. (3) \textbf{121} (2020), 51--82.

\bibitem{NT} N.~Nguyen, S.~Trivedi, \emph{Invariants and classification of simple function germs with respect to Lipschitz $\mathcal A$-equivalence}, J. Singul. \textbf{25} (2022), 348--360.

\bibitem{PT} F.~Pham, B.~Teissier, \emph{Fractions lipschitziennes d'une algèbre analytique complexe et saturation de Zariski}, prépublication, Centre de Mathématiques de l'École Polytechnique, 1969 (hal-00384928). English translation in: \emph{Introduction to Lipschitz Geometry of Singularities}, Lecture Notes in Math. 2280, Springer, 2020.

\bibitem{RV} M.~A.~S. Ruas, G.~Valette, \emph{$C^0$ and bi-Lipschitz $\mathcal K$-equivalence of mappings}, Math. Z. \textbf{269} (2011), 293--308.

\bibitem{S} J.~E. Sampaio, \emph{Bi-Lipschitz homeomorphic subanalytic sets have bi-Lipschitz homeomorphic tangent cones}, Selecta Math. (N.S.) \textbf{22} (2016), 553--559.

\bibitem{Wall} C.~T.~C. Wall, \emph{Finite determinacy of smooth map-germs}, Bull. London Math. Soc. \textbf{13} (1981), 481--539.
\bibitem{Wallbook} C.~T.~C. Wall, \emph{Singular Points of Plane Curves}, London Math. Soc. Student Texts 63, Cambridge University Press, 2004.
\end{thebibliography}
\end{document}